\documentclass{degm}

\newtheorem{theorem}{Theorem}[section]
\newtheorem{corollary}[theorem]{Corollary}

\newtheorem{lemma}[theorem]{Lemma}

\theorembodyfont{\rmfamily}

\newtheorem{remark}[theorem]{Remark}

\numberwithin{equation}{section}

\def\P{{\Bbb P}}
\def\H{{\Bbb H}}

\contact[conte@dm.unito.it]{Dipartimento di Matematica, Universit\`a di Torino,
Via Carlo Alberto 10, 10123 Torino, Italy}

\contact[marchisio@dm.unito.it]{Dipartimento di Matematica, Universit\`a di Torino,
Via Carlo Alberto 10, 10123 Torino, Italy}

\contact[murre@math.leidenuniv.nl]{Department of Mathematics, University of Leiden, P.O. Box 9512, 2300 RA Leiden, The Netherlands}

\begin{document}

\title{On Unirationality of Double Covers of Fixed Degree and Large Dimension; a Method of Ciliberto}

\author{Alberto Conte, Marina Marchisio and Jacob P. Murre}
\maketitle

\begin{flushright}
{{\Huge }\it To the memory of Paolo Francia}
\end{flushright}

\bigskip

\begin{abstract}
Following an idea of Ciliberto we show that double covers of projective $r$-space branched over an hypersurface of degree $2d$ are unirational provided $r$ is sufficiently big with respect to $d$.
\end{abstract}


\section{Introduction}
The notion of {\it unirationality} plays an important role in classical algebraic geometry in the works of M. Noether, Enriques and especially by Fano; see for istance Chap. IV of the book by Roth \cite{Ro}. At present the concept of {\it rationally connected} variety seems to become more and more important (see the recent book of Kollar \cite{Ko}), clearly unirational varieties are rationally connected but whether the latter concept is more general than the former is not yet known (ibid., problem 55). Irrespectively of the answer to this question, it remains an interesting geometrical problem to decide whether certain types of varieties are unirational (or not!). 

One the most striking results in this subject is a {\it theorem of U. Morin} from 1940 \cite{Mo1} saying that (always in characteristic zero) if $V=V_{r-1}(d) \subset \P^r$ is a {\it hypersurface} of degree $d$ in a projective $r$-space then there exists a constant $c(d)$ such that if $r\geq c(d)$ and if $V$ is ''sufficiently general'' then $V$ is unirational (see theorem \ref{Thm2} below for the precise statement). This theorem has been generalized to {\it complete intersections} by {\it Predonzan} in \cite{Pr2}. ''Modern'' treatments of the results of Morin and Predonzan were given in the papers of Ciliberto \cite{Ci}, Ramero \cite{Ra}, Paranjape-Srinivas \cite{PS} and we refer also to Chap.  10 in the book of Iskovskikh \cite{Isk}. Recently the result of Morin has been improved by Harris, Mazur and Pandharipande \cite{HPM} in the sense that ''sufficiently general'' has been relaxed to ''smooth'', this is an important improvement but - as far as we know - this result of \cite{HPM} for hypersurfaces has not yet been extended to complete intersections.

A natural question is whether the result of Morin can be extended to {\it double covers} $\pi:  W=W_r[2d,B]\longrightarrow \P^r$ of $\P^r$ ramified over an hypersurface $B=B_{r-1}(2d)\subset \P^r$ of degree $2d$; i.e., whether there exists a constant $\rho(d)$ such that if $r\geq \rho(d)$ and $B$ is ''sufficiently general'' the variety $W=W_r[2d,B]$ is unirational. In the above quoted paper \cite{Ci} Ciliberto has given a beautiful idea (osservazione 3.6) how to proceed to prove such a theorem (reducing it to a general criterion given by Morin in his Torino lecture of 1954 \cite{Mo2}). However, as Ciliberto remarks himself, the details of his outline depend upon a number of rather subtle verifications of ''algebraic'' nature. The purpose of this paper is to give these details and to prove the theorem for double covers; the precise statement is theorem \ref{Thm3} below. Our starting point is the theorem of Morin-Predonzan in the version of Ciliberto, see for the precise statement theorem \ref{Thm2} below.

There are in the theorems of Morin and Predonzan (at least) two important, but delicate, technical points: {\it firstly}, in order to specify the field over which the unirationality occurs one needs a ''sufficiently large'' linear space $L \simeq \P^q$ contained in $V$ and {\it secondly} the pair $(V,L)$ must be ''sufficiently general''. In the paper \cite{PS} and in the book of Iskovskikh one introduces the notion of ''general pair''. However it is difficult to control this notion of ''general pair''. On the other hand for the application of the results of Morin-Predonzan to the  case of double covers these technical aspects play an important (and in fact crucial) role. Therefore we have preferred to work with the precise notion of ''generic'' in the sense of Weil \cite{We} or Grothendieck \cite{EGA} (although technically differently framed the notions ''generic'' of Weil and Grothendieck are - of course - essentially the same). Working with ''generic'' it is important to distinguish between: ''$V$ generic over $K$'' (see subsection \ref{subsec21}) and in case $V$ contains a linear space $L$ ''$V$ generic over $K$ subject to containing $L$'' (see subsection \ref{subsec22}). We have given the precise definitions of these notions in section 2.

\section{Definition and preliminaries}\label{sec2}

\subsection{} \label{subsec21}

Let $K$ be a field of characteristic zero. Let $V = V_n$ be an irreducible variety defined over $K$ of dimension $n$. We recall that $V$ is called {\it unirational} if there exists a rational {\it dominant} map $f:   \P^n \longrightarrow V$ where $ \P^n$ is projective $n-$space, $f$ is defined over the algebraic closure $\overline{K}$ of $K$ (or better: usually over a finite extension $K'$ of $K$) and dominant means that the Zariski closure of the image $f( \P^n)$ in $V$ is $V$ itself. $V$ is {\it unirational over $K$} if moreover $f$ itself is also defined over $K$.

Let now $V=V_{r-1}(d) \subset  \P^r$ be a {\it hypersurface} of degree $d$, i.e. $V$ is defined in  $\P^r$ by an equation $F(X_0,\dots,X_r)=0$ homogeneous of degree $d$. We shall say that $V$ is {\it generic over $K$} if the coefficients of $F$ are {\it independent transcendental over $K$},  i.e. in the parameter space $\H(r,d)= \P^N$ with $N=\left(\begin{array}{cc}r+d\\d\\\end{array}\right)-1$ the $V$ corresponds to a generic point $\vartheta(V)$ over $K$ (in the terminology of Weil \cite{We}). Similarly if  $V=V_{r-m}(\underline{d}) \subset  \P^r$ is a {\it complete intersection} of multidegree $\underline{d}= (d_1,d_2,\dots,d_m) $ with $0<d_1\leq d_2 \leq \dots\leq d_m$ defined by equations $F_j(X_0,\dots,X_r)=0$ homogeneous of degree $d_j$ $(j=1,\dots,m)$ then we shall say that $V$ is {\it generic over $K$} if the coefficients of the $F_j$ are (mutually) independent transcendental over $K$ (note that such $V(\underline{d})$ are parametrized by a product of projective spaces $\H=\H(r,\underline{d})$  $= \P^{N_1}\times \dots \times \P^{N_m}$ with $N_j=\left(\begin{array}{cc}r+d_j\\d_j\\\end{array}\right)-1$).

\subsection{Linear Spaces contained in $V=V_{r-m}(\underline{d})$ }\label{subsec22} 

Let the field $K$ and $r$ and $\underline{d}$ be as above, and let $q$ be an integer such that $0<q<r$. We are interested in linear spaces $L = \P^{q}\subset \P^{r}$ which are contained in $V=V_{r-m}(\underline{d})$.  There is the following theorem of Predonzan \cite{Pr1} (see also \cite{Ci}, thm 2.1).

\begin{theorem}
\label{Thm1}
If $r$ and $\underline{d}= (d_1,d_2,\dots,d_m)$ are as above and $\underline{d}\neq (1,1,\dots,1,2)$ then $V=V_{r-m}(\underline{d})$ contains a linear space $L=\P^{q}$ if

\begin{equation}
(r-q)(q+1)\geq \sum_{i}\left(\begin{array}{cc}d_j+q\\q\\\end{array}\right).
\label{1}
\end{equation}
\end{theorem}

Consider the incidence correspondence
\bigskip
\[
\begin{array}{cccc}I=&I(r,\underline{d},q)&\stackrel{p_2}{\longrightarrow}&\H(r,\underline{d})\\ [.11in]
\empty & \mbox{\scriptsize{\it {p}$_1$}}\downarrow &\empty&\empty\\[.11in]
\empty &{\bf Gr}(q,r)&\empty&\empty\\
\end{array}
\]
\bigskip

\noindent where, as usual, ${\bf Gr}(q,r)$ is the Grassmannian of the $\P^{q}\subset \P^{r}$, $\H(r,\underline{d})$ is the variety parametrizing the complete intersections of multidegree $\underline{d}$ in $\P^r$ and

\[
I=I(r,\underline{d},q)=\{(V,L); L\in {\bf Gr}(q,r), V \in \H(r,\underline{d})\,\, \mbox{ and } \,\, L\subset V\}.
\]

\bigskip

Given $V$, put $F(V,q)=p_2^{-1}(V)=\{L \in {\bf Gr}(q,r); L\subset V\}$; this is the {\it Fano variety} of  $\P^q$'s in $V$.

It is well-known (and in fact very easy to see) that $I$ is irreducible over $K$ and of dimension

\[
(r-q)(q+1)+ \sum_{i}\left(\begin{array}{cc}d_i+r\\r\\\end{array}\right)- \sum_{i}\left(\begin{array}{cc}d_i+q\\q\\\end{array}\right)-m
\]

\noindent and $p_2$ is {\it onto} only if, and by the above theorem \ref{Thm1} of Predonzan in fact if, the inequality (\ref{1}) from above holds.

Given the field $K$, let $(V,L) \in I$. We shall say that {\it the pair} $(V,L)$ is {\it generic over $K$} if $(V,L)$ is a generic point of $I$. From the above we have immediately

\begin{lemma}
Assume $($\ref{1}$)$ holds. Let $(V,L) \in I$. Then the following are equivalent

\noindent $1.$ $(V,L)$ is a {\it generic pair} over $K$

\noindent $2.$ $L \in {\bf Gr}(q,r)$ is generic over $K$ and $V$ is generic over $K(l(L))$ in the fibre $p_1^{-1}(L)$

\noindent $3.$ $V\in \H(r,\underline{d})$ is generic over $K$, $F(V,q)$ is irreducible over $K(\vartheta(V))$ and $L$ is generic over $K(\vartheta(V))$ in $F(V,q)$.
\end{lemma}

\begin{remark}
By $K(l(L))$, resp. $K(\vartheta(V))$, we denote the field of definition over $K$ of $L$, resp. of $V$ (i.e., obtained by adjoining to $K$ the ratios of the pl\"ucker coordinates, resp. the ratios of the coefficients of the equations).
\end{remark}

Given the field $K$, let $L\in {\bf Gr}(q,r)$; let $(V,L)\in I$. We shall say that $V$ is {\it generic over $K$ subject to containing $L$} if $V$ is generic over $K(l(L))$ in the fibre $p_1^{-1}(L).$ Note that this fibre itself is (isomorphic to) a product of projective spaces of type $\P^{M_1}\times \dots \times \P^{M_m}$ with $M_i = \left(\begin{array}{cc}d_i+r\\r\\\end{array}\right)- \left(\begin{array}{cc}d_i+q\\q\\\end{array}\right)-1$ (of course in general $V$ will no longer by generic over $K$, i.e., no longer generic in $\H(r,\underline{d})$).

Given such $L=\P^q\subset \P^r$, say $L=L_0$, we can choose homogeneous coordinates $(Z_0,Z_1,\dots,Z_q,Y_{q+1},\dots,Y_r)$ in $\P^r$ such that $L_0$ is given by

\begin{equation}
Y_{q+1}=Y_{q+2}=\dots = Y_r=0
\label{2}
\end{equation}

\bigskip
\noindent and we can use then $(Z_0,\dots,Z_q)$ as homogeneous coordinates in $L_0$.

\bigskip

\noindent {\bf Note.} If $L_0$ is defined over $K$ then we can make this projective coordinate transformation over $K$ itself, however if $L_0$ {\it should not be defined over $K$ } (which of course is the general situation) we need to make base extension $K(l(L_0))\supset K$ and perform this transformation over $K(l(L_0))$.

\bigskip

If we {choose coordinates in this way} then the equations of $V=V_{r-m}(\underline{d}) \subset \P^r$ take the following shape:

\begin{equation}
{\cal{G}}_j(Z_0,\dots,Z_q,Y_{q+1},\dots,Y_r)=\sum_{\begin{array}{cc}I\\0\leq|I|<d_j\\\end{array}}Z^I{\cal{G}}_{j,I}(Y)=0 \,\,\,\,(j=1,\dots,m),
\label{3}
\end{equation}

\noindent where $I=(i_0,i_1,\dots,i_q)$ and $Z^I=Z_0^{i_0}Z_1^{i_1}\dots Z_q^{i_q},$ $I=\sum^q_{\rho=0}i_{\rho}$, i.e., the total degree in the $Z_i$'s and ${\cal{G}}_{j,I}(Y)$ is homogeneous in the  $Y_i$'s of degree $d_j-|I|$. Note that the condition that the $Z^I$ with $|I|=d_j$ do {\it not} occur is precisely the condition that $L_0\subset V$. Also note that the coefficients of the ${\cal{G}}_j(Z,W)$ are in $K(l(L_0))$.

Now by counting the dimensions we have immediately the following lemma

\begin{lemma}
Given $K$ and $L_0 \in {\bf Gr}(q,r),$  let $(V,L_0) \in I(r,\underline{d},q)$. Then  the following are equivalent:

\noindent $1.$ $V$ is generic over $K$ subject to containing $L_0$,

\noindent $2.$ in the above equations $($\ref{3}$)$ the (ratios of the) coefficients are independent transcendental over $K(l(L_0))$.
\end{lemma}\label{lem2}

\begin{remark}
Of course we mean not only for one separate $j$, but for all the $j$'s together.
\end{remark}

\section{Theorems of Morin and Predonzan, in version of Ciliberto}\label{sec3}

\begin{theorem}  {\rm (Morin, Predonzan, Ciliberto)}
\label{Thm2}

\noindent Given a field $K$ of characteristic zero and $\underline{d}=(d_1,d_2,\dots,d_m)$ with $0<d_1\leq  d_2 \leq \dots \leq d_m$. Then there exist integers $c(\underline{d})$ and $q(\underline{d})$ such that if $r\geq c(\underline{d})$ and $L_0 \in {\bf Gr}(q(\underline{d}),r)$ and if $V=V_{r-m}(\underline{d}) \subset \P^r$ is  generic over $K$ subject to containing $L_0$ then $V$ is unirational over the field $K(l(L_0),\vartheta(V))$.
\end{theorem}

For $\underline{d}=d$ this is the theorem of Morin, for general $\underline{d}$ we get the theorem of Predonzan. The above theorem is Corollary 2.5 of Ciliberto \cite{Ci}, his notion ''generica su $K(l(L_0))$'' is indeed precisely what we call ''$V=V_{r-m}(\underline{d})\subset \P^r$ is generic over $K$ subject to containing $L_0$'' as we see from his remark on page 180 of his paper on the equations (2.7) in his paper. Finally in Ciliberto's notation $c(\underline{d})=s_n$ and $q(\underline{d})=s_{n-1}$ (see his corollary 2.5); however remark that $n$ in the paper of Ciliberto is determined by $\underline{d}$, and $s_n$ and $s_{n-1}$ are determined by $\underline{d}$, see his definitions on page 177.

\bigskip

We shall use the theorem of Morin and of Predonzan in the above precise version of Ciliberto.

\section{Double covers of $\P^r$ and unirationality}

\subsection{}

Let, as above, $K$ be a field of characteristic zero. Let $B=B_{r-1}(2d)\subset\P^r$ be a hypersurface of degree $2d$ in $\P^r$ and let $\pi:  W= W_r[2d,B]\longrightarrow \P^r$  be the double cover of $\P^r$ branched over $B$. We want to prove the following:

\begin{theorem}\label{Thm3} 
Given  $K$ and $d>2$, there exists a constant $\rho(d)$ such that if $r\geq \rho(d)$ and if $B_{r-1}(2d)\subset  \P^r$ is generic over $K$ then the double cover $W[2d,B]$ of $ \P^r$ is unirational, and in fact unirational over a finite extension $K^{\#}$ of $K(\vartheta(B))$.
\end{theorem}

\section{Proof of the theorem}

\subsection{Preparations and beginning of proof}

Let $d\geq 3$ and put $\underline{d}=(1,2,\dots,2d-2)$. Let furthermore $q(\underline{d})$ and $c(\underline{d})$ be the integers occuring in theorem \ref{Thm2} (i.e., the integers $s_{n-1}$ and $s_n$ in Ciliberto's corollary 2.5 of \cite{Ci}).

\begin{lemma}
\label{lem3}
For $d \geq 3$ the $q(\underline{d})\geq 2d-2.$
\end{lemma}

\begin{proof}
Elementary, but cumbersome. We leave it to the reader but we make some remarks. The proof goes by induction and starting with $d=3$. Examining the expressions (2.3) on page 177 of Ciliberto's paper we see that in our case his $n=2d-2$. For $d=3$ we get $q(3)=25$ and since $25>4$ we can start. Next for the induction step, passing from $(d-1)$ to $d$ we see that the $q(\underline{d})$ increases by at least $2$.
\end{proof}

We take 

\begin{equation}
\label{51}
q= q(\underline{d})+1=:\rho''(d).
\end{equation}

\bigskip
Now we {\it introduce a constant $\rho_1(d)$} as follows: by theorem \ref{1} there exists a constant $c^*(2d,q)$ such that if $V=V_{n-1}(2d)\subset \P^n$ then $V$ contains a linear space $L\simeq \P^q$ provided $n\geq c^*(2d,q)$. Now take

\begin{equation}
\label{52}
\rho_1(d)=c^*(2d,\rho''(d)).
\end{equation}

\begin{lemma}
\label{lem4}
Given the field $K$ and $d\geq 3$, take $q$ from $($\ref{51}$)$. Let $r\geq \rho_1(d)$. Let $B=B_{r-1}(2d)\subset \P^r$ be generic over $K$. Then there exists a linear space $L_0 \simeq \P^q$ such that

(i) $L_0 \subset B$ and $L_0$ is defined over a finite extension $K^{\#}$ of $K(\vartheta(B))$

(ii) $B$ is generic over $K$ subject to containing $L_0$ (recall that this notion is defined in subsection \ref{subsec22} and it means that $B$ is generic over $K(l(L_0))$ in the fibre $p^{-1}_1(L_0)$ of the diagram below).
\end{lemma}

\begin{proof}
Consider the incidence diagram (like in subsection \ref{subsec22})

\[
\begin{array}{cccc}I=&I(r,2d,q)&\stackrel{p_2}{\longrightarrow}&\H(r,2d)\ni B\\[.11in] \empty &\mbox{\scriptsize{\it {p}$_1$}}\downarrow &\empty&\empty\\[.11in]\empty &{\bf Gr}(q,r) \supset F(B,q)&\empty&\empty\\
\end{array}
\]

\bigskip
\noindent Consider the Fano variety $F(B,q) = p_2^{-1}(B)=\{L\in {\bf Gr}(q,r);L\subset B\},$ this variety is non-empty since $r\geq \rho_1(d)$, it is defined over $K(\vartheta(B))$ and as we have remarked in subsection \ref{subsec22} it is irreducible over $K(\vartheta(B))$ and of dimension

\[
\dim F(B,q)=(q+1)(r-q)-\left(\begin{array}{cc}q+2d\\2d\\\end{array}\right)=: b(r,d).
\]

\bigskip
Now intersect this Fano variety with the linear section (via hyperplanes in the Pl\"ucker embedding) ${\cal{L}}$ of ${\bf Gr}(q,r)$ of codimension $b(r,d)$, i.e., of dimension $\left(\begin{array}{cc}q+2d\\2d\\\end{array}\right)$, {\it defined over $K$} and sufficiently general such that the intersection consists of points (this is possible, $K$ has infinitely many elements hence we can choose the hyperplane ''sufficiently general''). These points are defined over the algebraic closure $\overline{K(\vartheta(B))}$ of $K(\vartheta(B))$; take one of them, say, $L_0$ so $L_0 \in F(B,q) \cap {\cal{L}}$. Now we claim that $B$ is (still) generic in $p_1^{-1}(L_0)$. In fact we have

\[
\left(\begin{array}{cc}r+2d\\2d\\\end{array}\right)-1=\mbox{ trdeg}_K K(\vartheta(B))=
\mbox{ trdeg}_K K(\vartheta(B),l(L_0))=
\]
\[
\mbox{ trdeg}_K K(l(L_0))+\mbox{ trdeg}_{K(l(L_0))}K(l(L_0),\vartheta(B))\leq\left(\begin{array}{cc}q+2d\\2d\\\end{array}\right)+\dim p_1^{-1}(L_0)=
\]
\[
\left(\begin{array}{cc}q+2d\\2d\\\end{array}\right) + \left\{\left(\begin{array}{cc}r+2d\\2d\\\end{array}\right) -1 -\left(\begin{array}{cc}q+2d\\2d\\\end{array}\right)\right\}.
\]

\bigskip
\noindent (Note: we use $L_0 \in {\cal{L}},$ $ \dim {\cal{L}} =\left(\begin{array}{cc}q+2d\\2d\\\end{array}\right)$). Hence we must have equality, in particular $\mbox{ trdeg}_{K(l(L_0))} K(l(L_0),$ $\vartheta(B)) = \dim p_1^{-1}(L_0)$, i.e., $B$ is generic over $K(l(L_0))$ in the fibre $p_1^{-1}(L_0)$.
This completes the proof of the lemma.
\end{proof}

\subsection{Continuation. Choice of $\rho(d)$}

Take

\begin{equation}
\label{53}
\rho(d):=\max\{c(\underline{d}), \rho_1(d)\}+1,
\end{equation}

\bigskip
\noindent where $c(\underline{d})$ is the constant in theorem \ref{Thm2} with $\underline{d}=(1,2,\dots,2d-2)$ and $\rho_1(d)$ is from (\ref{52}).

Now let $r\geq \rho(d)$, take $B=B_{r-1}(2d)\subset \P^r$ {\it generic over $K$} and let $W_r[2d,B]$ be a {\it double cover} of $\P^r$ {\it branched} over $B$. {\it We must prove that } $W_r[2d,B]$ is {\it unirational} and in fact unirational over a finite extension of $K(\vartheta(B))$.

Let $q$ be the integer from (\ref{51}). Since $r\geq \rho_1 (d)$ we can, by lemma \ref{lem4}, find a linear space $L_0\subset B$ of dimension $q$ satisfying the conditions of lemma \ref{lem4}, {\it in particular $B$ is generic over $K$ subject to containing $L_0$} (also: $L_0$ is defined over a finite extension of $K(\vartheta(B))$).

Fix moreover in $L_0$ a linear space $M_0\subset L_0$ of dimension $(2d-2)$ (possible since $q\geq  (2d-2)$ by lemma \ref{lem3}) and take $M_0$ to be defined over $K(l(L_0))$.

After we have choosen such $L_0\subset B=B_{r-1}(2d)\subset \P^r$ we make a projective coordinate transformation, defined over $K(l(L_0))$, such that we have homogeneous coordinates $Z_0,\dots , Z_q ,Y_{q+1},\dots,Y_r$ such that $L_0$ is given by (see subsection \ref{subsec22})

\begin{equation}\label{54}
Y_{q+1}=\dots=Y_r=0
\end{equation}

\noindent and we use $Z_0, Z_1, \dots, Z_q$ as homogeneous coordinates in $L_0$.

We can take for $M_0\subset L_0$ the space defined by

\begin{equation}\label{55}
Z_{2d-1}=\dots=Z_q=Y_{q+1}=\dots=Y_r=0
\end{equation}

The equation of $B_{r-1}=B_{r-1}(2d)\subset \P^r$ has now the following form (see in section \ref{sec2} the equation (\ref{3}))

\begin{equation}\label{56}
\mathcal G(Z,Y)=\sum_{\begin{array}{cc}I\\0\leq|I|<2d\\\end{array}} Z^I \mathcal G_I(Y)=0
\end{equation}

\noindent where $I=(i_0,\dots,i_q)$,  $|I|=\sum_{\sigma=0}^q i_{\sigma}$ and $\mathcal G_I(Y)$ is homogeneous in the $Y_j$'s  of degree $2d-|I|$.
The coefficients of $\mathcal G_{i_0\dots i_q}(Y)$ are $b_{i_0 \dots i_q j_{q+1}\dots j_{\sigma}}$ (the sum of all the indices is 2d).
Since $B$ is generic over $K$ subject to containing $L_0$ we have, by lemma \ref{lem2} that the

\begin{equation}
b_{i_o\dots i_q j_{q+1}\dots j_\sigma} \textrm{are independent transcendentals over }K(l(L_0)).
\end{equation}\label{57}

\noindent Next we choose a hyperplane $H_0$ in $\P^r$ not containing $M_0$ and, a fortiori, therefore not containing $L_0$.
Let this hyperplane be $\{Z_0=0\}$.
Put $L_0^{*}=L_0\cap H_0$ and $M_0^{*}=M_0\cap H_0$.
So the equations of $L_0^*$ are 

\begin{equation}\label{58}
Z_0=0,Y_{q+1}=Y_{q+2}=\dots=Y_r=0
\end{equation}

\noindent and similar for the $M_0^{*}$.

\subsection{The varieties $F_\eta$ and $F_{\eta}^*$}

Take in $M_0$ a point $\eta$ generic over $K(l(L_0),\vartheta(B))$, so $\eta$ has coordinates (see the equations (\ref{55}))

\begin{equation}\label{59}
\eta =(\eta_0=1,\eta_1,\dots,\eta_{2d-2},0,\dots ,0,\dots ,0)
\end{equation}

\noindent with $\eta_1,\dots,\eta_{2d-2}$ independent transcendentals over $K(l(L_0),\vartheta(B))$.
For the sake of simplicity of notation we sometimes write

\begin{equation}
\eta =(\eta_0=1,\eta_1,\dots,\eta_{2d-2},\eta_{2d-1}=0,\dots,\eta_q=0,0,\dots ,0).
\end{equation}\label{59'}

Now consider the (higher) {\it polar varieties} of $\eta$ with respect to $B$, up to the $(2d-2)$-th one, i.e., the varieties

\begin{equation}
\Delta_{\eta}^1(B)=0,\Delta_{\eta}^2(B)=0,\dots ,\Delta_{\eta}^{2d-2}(B)=0,
\label{510}
\end{equation}

\noindent where $\Delta_{\eta}^{j}(B)$ is defined by the equation

\begin{equation}
\Delta_\eta^j(B)=\left.\left(Z_0\frac{\partial}{\partial Z_0}+\dots +Z_q\frac{\partial}{\partial Z_q}+Y_{q+1}\frac{\partial}{\partial Y_{q+1}}+\dots +Y_r\frac{\partial}{\partial Y_r}\right)^{(j)}{\cal{G}}(Z,Y)\right|_\eta=0.
\label{511}
\end{equation}

So in particular $\Delta_{\eta}^1(B)=0$ is the tangent space to $B$ in $\eta$, $\Delta_{\eta}^1(B)\cap\Delta_{\eta}^2(B)$ is the tangent cone to $B$ in $\eta$,etc.

Let

\begin{equation}
F_\eta=\Delta_\eta^1(B)\cap\Delta_\eta^2(B)\cap\dots\cap\Delta_\eta^{2d-2}(B).
\label{512}
\end {equation}

So $F_{\eta}$ consists of the lines $m$ through $\eta$ such that 

\begin{equation}
m\cap B=(2d-1)\eta+\beta
\label{513}
\end{equation}

\noindent with $\beta\in B\cap F_{\eta}:=B_{\eta}^*$.
Also put $F_{\eta}^*=F_{\eta}\cap H_0$ and clearly $F_{\eta}^*$ is the projection of $F_{\eta}^*$ from the point $\eta$ and is birationally equivalent with $F_{\eta}^*$ over the field $K(l(L_0),\vartheta(B),\eta)$. Also clearly $F_{\eta}^*$ is a variety of type $V_{r-2d+1}(1,2,\dots,2d-2) \subset H_0$ and it contains the $(q-1)$-dimensional linear space $L_0^* \,\,\,(= L_0\cap H_0)$. 

The following lemma is crucial (but subtle!).

\begin{lemma}\label{lem5}
$F_{\eta}^*$ is a variety of type $V_{r-2d+1}(1,2,\dots,2d-2)$ in $H_0$ and $F_{\eta}^*$ is generic over $K(l(L_0),\eta)$ subject to containing $L_0^*$ (in the sense of subsection \ref{subsec22}).
\end{lemma}

\bigskip
\noindent {\bf Note.} This lemma plays in our case a role analogous to the assertion i) on page 186 of Ciliberto's paper \cite{Ci}. Since the lemma is crucial for our purpose we give a full proof (with all details).

\begin{proof}
We start with two remarks.

\begin{remark}\label{rem1}
Working first (for simplicity) over the field $K':= K(l(L_0))$, we can consider {\it over this field} a projective coordinate transformation $T$ in the linear space $M_0$ of the type $\stackrel{\sim}{Z}_0=Z_0,\stackrel{\sim}{Z}_l=\sum_{i=0}^{2d-2}a_{il}Z_i=T(Z_0,\dots,Z_{2d-2})$ with $a_{il} \in K'$. Of course we can also consider - if we prefer - the $T$ as a linear coordinate transformation in $\P^r \supset L_0 \supset M_0$ itself (leaving the other coordinates unchanged). (Note also that it leaves the equation $Z_0=0$ of $H_0$ unchanged.)

$T$ operates now on the parameter space $\H(r,2d)$ of the hypersurfaces of degree $2d$ in $\P^r$ (see the diagram in the proof of lemma \ref{lem4}), leaving the space $p_1^{-1}(L_0)$ fixed. It transforms the point $B\in p_1^{-1}(L_0)$ into a point $B^T \in p_1^{-1}(L_0)$ corresponding to the (new) equation obtained from (\ref{56}):

\begin{equation}
{\mathcal G}^T(\stackrel{\sim}{Z},Y)={\mathcal G}(T^{-1}(\stackrel{\sim}{Z}),Y) =\sum_{\begin{array}{cc}I\\|I|<2d\\\end{array}} T^{-1}(\stackrel{\sim}{Z})^I{\mathcal G}_I(Y)=\sum_{\begin{array}{cc}I\\|I|<2d\\\end{array}} \stackrel{\sim}{Z}^I{\mathcal G}_I^T(Y).
\end{equation}\label{56uno}

\bigskip

\noindent{\bf Claim}. {\it The $B^T$ is (still) generic over $K$ subject to containing $L_0$ (in the sense of section \ref{sec2}, i.e., $B^T \in p_1^{-1}(L_0)$ is a generic point of $p_1^{-1}(L_0)$ over $K(l(L_0))).$}

\bigskip
\noindent {\it Proof of the Claim.} The coefficients $b^T_{i_0i_1\dots i_qj_{q+1}\dots j_T}$ of ${\mathcal G}^T(\stackrel{\sim}{Z},Y)$ are linear over $K'=K(l(L_0))$ in the coefficients of ${\mathcal G}^T(Z,Y)$, hence $K'(\vartheta(B^T))\subset K'(\vartheta(B))$. Applying the inverse of $T$ we get $K'(\vartheta(B))=K'(\vartheta(B^T)).$
\bigskip
\end{remark}

\begin{remark}\label{rem2}
Let $\eta \in M_0$ be a generic point of $M_0$ over $K(l(L_0),\vartheta(B))$, then ''conversely'' $B$ is still generic over the field $K(\eta)$ subject to containing $L_0$ (count transcendence degrees).

Hence if we replace the field $K'=K(l(L_0))$ by the field $K'':=K(\eta, l(L_0))=K'(\eta)$ then we can apply the above remark \ref{rem1} to $B$ and the field $K''.$
\end{remark}

After the two remarks we proceed as follows.

Apply over the field $K''=K(\eta, l(L_0))$ in $M_0$ the projective coordinate transformation $T$ such that $\stackrel{\sim}{Z}_0=Z_0$ and such that the point $\eta=(1,\eta_1,\dots,\eta_{2d-2},0,\dots,0)$ becomes the point $\eta_0=(1,0,\dots,0)$.

Then the hyperplane $B=B_{r-1}(2d)\subset \P^r$ corresponds in the parameter space $\H(r,2d)$ with the point $B^T$ corresponding to the equation

 \begin{equation}
{\mathcal G}^T(\stackrel{\sim}{Z},Y) =\sum_{\begin{array}{cc}I\\|I|<2d\\\end{array}} \stackrel{\sim}{Z}^I{\mathcal G}_I^T(Y)=\sum_{\begin{array}{cc}I\\|I|<2d\\I+J=2d\\\end{array}} \stackrel{\sim}{b}_{I,J}\stackrel{\sim}{Z}^IY^J
\end{equation}\label{56unouno}

\noindent and by remark \ref{rem1} the coefficients $\stackrel{\sim}{b}_{I,J} (\in K''(\vartheta(B)))$ are all {\it independent transcendental} over $K''=K(\eta,l(L_0))$ (and moreover also note that in {\it the allowed range}, i.e., $0\leq|I|<2d,$ {\it all coefficients occur}).

Now we are going to rearrange them according to {\it decreasing} powers of $\stackrel{\sim}{Z}_0=Z_0.$ (Note since $\eta_0\in B \,\,\,\,Z_0^{2d}$ does not occur.)

We get

 \begin{equation}\label{56due}
\begin{array}{r}
{\mathcal G}^T(\stackrel{\sim}{Z},Y) =Z_0^{2d-1}\Phi_1(Y_{q+1},\dots,Y_r)+Z_0^{2d-2} \Phi_2(\stackrel{\sim}{Z}_1,\dots,\stackrel{\sim}{Z}_q,Y_{q+1},\dots,Y_r)+\\
 \dots +Z_0^{2d-s} \Phi_s(\stackrel{\sim}{Z},Y)+\dots\\
\end{array}
\end{equation}

\noindent where the $\Phi_s(\stackrel{\sim}{Z},Y)$ are homogeneous polynomials in the $\stackrel{\sim}{Z}_1,\dots,\stackrel{\sim}{Z}_q,Y_{q+1},\dots,Y_r$ of degree $s$, {\it subject to the condition that they don't contain terms in the $\stackrel{\sim}{Z}_1,\dots,\stackrel{\sim}{Z}_q$ alone} (this $\Longleftrightarrow |I| <2d \Longleftrightarrow B$ contains $L_0$) (in particular $\Phi_1(\stackrel{\sim}{Z},Y)=\Phi_1(Y)$). The coefficients in the above equation are the $\stackrel{\sim}{b}_{I,J}$ from above, hence by what we did say above they are independent transcendental over $K''=K(l(L_0),\eta)$ and all of them occur in the allowed range.

However now (and compare with Ciliberto page 187) the $s$-th  polar of the point $\eta_0=(1,0,\dots,0)$ with respect to $B$ is precisely given by

\[
\Delta^s_{\eta_0}(B^T)=\Phi_s(\stackrel{\sim}{Z}_1,\dots,\stackrel{\sim}{Z}_q,Y_{q+1},\dots,Y_r)=0
\]

\noindent and the $F^*_{\eta}=(F^T)^*_{\eta_0}$ is given by $\Phi_1=\Phi_2=\dots=\Phi_{2d-2}=0$.

Hence by what we have said above about the coefficients $ \stackrel{\sim}{b}_{I,J}$ in (\ref{56due}) it follows that the $F^*_{\eta}$ is indeed a $V(1,2,\dots,2d-2)\subset H_0$ generic over $K''=K(\eta,l(L_0))$ subject to containing $L^*_0$ ( this latter condition is precisely equivalent to the condition that there are in the $\Phi_s (s=1,\dots,2d-2)$ no terms in the $\stackrel{\sim}{Z}_i$'s alone).
\end{proof}

\begin{corollary}\label{cAA}
If $ r\,\geq\,\rho(d) $ $($see $($\ref{53}$))$ then $\,\,F^*_{\eta} $ is unirational over the field $ K(l(L_0),\eta ,\vartheta(B))$.
\end{corollary}

\begin{proof}
By lemma \ref{lem5} and since $r-1\geq c(\underline{d})$ we can apply the theorem of Predonzan-Ciliberto, i.e. theorem \ref{Thm2}.
\end{proof}

Next, in order to simplify notation, let us write

\begin{equation}\label{519}
K_1:=K(\vartheta(B))=K(\vartheta(W))\subset K'_1:= K(\vartheta(B),l(L_0)).
\end{equation}

\bigskip
From the corollary \ref{cAA} we have that the function field $K'_1(\eta)(F^*_{\eta})$ of $F^*_{\eta}$ over $K'_1(\eta)$ is contained in a purely transcendental extension. Therefore if $\xi \in F^*_{\eta}$ is generic over $K'_1(\eta)$ then we have

\begin{equation}\label{520}
K'_1(\eta,\xi)\cong K'_1(\eta)(F^*_{\eta})\subset K'_1(\eta,t_1,\dots,t_{r-2d+1})
\end{equation}

\noindent where the $t_1,\dots,t_{r-2d+1}$ are independent transcendental over $K'_1(\eta)$.

\subsection{Construction of a unirational family of lines in $\P^r$}\label{subsec54}

Let $S \subset M_0\times H_0$ be the Zarisky closure over $K'_1$ of the point $(\eta,\xi)$ where $\eta\in M_0$ is generic over $K'_1$ and $\xi \in F^*_{\eta}$ is generic over $K'_1(\eta)$. From (\ref{520}) it follows that $S$ is a variety which is {\it unirational} over $K'_1$, since we have

\begin{equation}\label{521}
K'_1(S):\cong K'_1(\eta,\xi)\subset K'_1(\eta_1,\dots,\eta_{2d-2},t_1,\dots,t_{r-2d+1})
\end{equation}

\noindent and the right hand side of (\ref{521}) is a purely transcendental extension of $K'_1$. Moreover, clearly $\dim S= r-1$.

Now we consider the line $m=\langle\eta,\xi\rangle$ spanned by $\eta$ and $\xi$. Then $\langle\eta,\xi\rangle \subset F_{\eta}$ and by (\ref{513}) we have 

\begin{equation}\label{522}
\langle\eta,\xi\rangle \cap B = (2d-1) \eta +\beta
\end{equation}

\noindent with $\beta \in B^*_{\eta} = B\cap F_{\eta}$.

\begin{lemma}\label{lem7}
When $\eta \in M_0$ is generic over $K'_1$ and $\xi \in F^*_{\eta}$ is generic over $K'_1(\eta)$, then $\beta\in B$ is generic over $K'_1$.
\end{lemma}

\begin{proof}
This is assertion ii) of page 186 of Ciliberto \cite{Ci}. For the convenience of the reader we repeat - in our language and notation - the argument here.

Let us denote during the proof of this lemma by $K'$ the field $K':=K(l(L_0)).$

With the notations and points from our lemma we have the following inclusions of fields

\[
\begin{array}{cccr}K'(\vartheta(B))&\longrightarrow&K'(\vartheta(B),\eta)\empty\\[.11in] 
\downarrow &\empty&\downarrow&\empty\\[.11in]
K'(\vartheta(B),\beta)&\longrightarrow&K'(\vartheta(B),\eta,\beta)&=K'(\vartheta(B),\eta,\xi)\\
\end{array}
\]

\bigskip
\noindent We have $\mbox{ trdeg}(K'(\vartheta(B),\eta):K'(\vartheta(B)))=2d-2$ and $\mbox{ trdeg}(K'(\vartheta(B),\eta,\xi):K'(\vartheta(B),\eta)))=r-2d+1,$ hence $\mbox{ trdeg}(K'(\vartheta(B),\eta,\beta):K'(\vartheta(B)))=r-1.$ Furthermore clearly $\mbox{ trdeg}(K'(\vartheta(B),\beta):K'(\vartheta(B))\leq r-1$ and hence in order {\it to prove the lemma} we have to prove that $K'(\vartheta(B),\eta,\beta)$ is an {\it algebraic extension} of $K'(\vartheta(B),\beta)$. 

The proof goes by contradiction. Suppose that 

\[
(*) \mbox{ trdeg}(K'(\vartheta(B),\eta,\beta):K'(\vartheta(B),\beta))>0.
\]

\noindent Now comes the nice {\it specialization argument} of Ciliberto! Consider the following variety

\[
J\subset p_1^{-1}(L_0)\times M_0\times \P^r \subset \H(r,2d)\times M_0\times \P^r
\]

\noindent with $J$ the Zariski locus over $K'$ of the point $(\vartheta(B),\eta ,\beta)$ and let $J_{13}= \mbox{ pr}_{13} J.$ By construction these are irreducible varieties over $K'$.

Now our assumption $(*)$ about the positive transcendence degree means that over the generic point $(\vartheta(B),\beta)$ of $J_{13}$ the fibre $p_{13}^{-1}(\vartheta(B),\beta)$ of $J$ has positive dimension. So then {\it the fibres over all points of $J_{13}$ have positive dimension.} Consider now the following special point $(\vartheta(B^*),\eta^*, \beta^*)$ of $J$ with (with the coordinates of (\ref{55})):\\
$B^*\subset \P^r$ is the hypersurface with equation

\[
Z^{2d-1}_{2d-3}Y_r+Z^{2}_{0}Y_r^{2d-2}+Z^{3}_{1}Y_r^{2d-3}+\dots+Z^{2d-2}_{2d-4}Y_r^{2}=0
\]

\noindent (note that $B^* \supset L_0$, in fact $B^* \supset (Y_r=0)$)

\[
\begin{array}{ccc}
\eta^*&=&(0,\dots,0,1,0,\dots,0)\\[.11in]
\empty&\empty& \,\,\mbox{place } \uparrow2d-2\\[.11in]
\beta^*&=&(0,\dots,0,1).\\
\end{array}
\]

\noindent Then one checks that indeed $\beta^* \in B^* \cap F_{\eta^*}$ (i.e. that $\beta^*$ is on the polars of $\eta^*$ with respect to $B^*$), and hence indeed $(\vartheta(B^*),\eta^*,\beta^*) \in J$ and hence $(\vartheta(B^*),\beta^*) \in p_{13}(J)=J_{13}.$ However now one also checks by direct computation (via the equations $\Delta^j_{\eta}(B))$ that $\eta^*$ is the {\it only} point in the fibre $p_{13}^{1}(\vartheta(B^*),\beta^*)$! (counted in fact $(2d-1)!$ times). However this contradicts our assumption $(*)$, which proves the lemma.
\end{proof}

\begin{corollary}\label{cor8}
Let $\eta$ and $\xi$ be as above and let $\alpha$ be a point on $\langle\eta,\xi\rangle$ generic over $K'_1(\eta,\xi)$. Then $\alpha$ is a generic point of $\P^r$ over $K'_1.$
\end{corollary}

\begin{proof}
By lemma \ref{lem7} the Zariski closure of $\alpha$ over $K'_1$ contains $B$; since it clearly contains $B$ {\it strictly} it must be $\P^r.$
\end{proof}

\subsection{Construction of a unirational family of rational curves on $W$. End of the proof}

Returning to the double cover $\pi: W_r[2d,B] \longrightarrow \P^r$ for $r\geq \rho(d)$ (from (\ref{53})), consider on $W=W_r[2d,B]$ the curve $C_{\eta\xi}=\pi^{-1}(\langle\eta,\xi\rangle)$ where $\eta$ and $\xi$ are as in corollary \ref{cor8}.

\begin{lemma}\label{lem9}
The curve $C_{\eta\xi}$ is a rational curve, rational over the field $K'_1(\eta,\xi)$.
\end{lemma}

\begin{proof}
The curve is a double cover of the line $\langle\eta,\xi\rangle$ branched only over the two points $\eta$ and $\beta$ (the $\beta$ from (\ref{522})). So it is a rational curve. Moreover it is rational over the field $K'_1(\xi,\eta)$ since it contains a rational point $\stackrel{\sim}{\eta}$ over this field where $\stackrel{\sim}{\eta}\in W$ such that $\pi(\stackrel{\sim}{\eta})=\eta$.
\end{proof}

Let $\omega\in C_{\eta\xi}$ be a generic point of this curve over $K'_1(\eta,\xi)$ and consider in $M_0\times H_0\times W$ the Zariski closure $W'$ over $K'_1$ of the point $(\eta,\xi,\omega)$. Then we have a diagram

\[
\begin{array}{cccc}
\empty&W'&\stackrel{p_3}{\longrightarrow}&W\\ [.11in]
\empty&\mbox{\scriptsize{\it {p}$_{12}$}}\downarrow &\empty&\empty\\[.11in]
\empty &S&\empty&\empty\\
\end{array}
\]

\bigskip
\noindent where $S$ is the variety from subsection \ref{subsec54}.

\begin{lemma}\label{lem10}
$W'$ is unirational over $K'_1$ and $p_3$ is onto.
\end{lemma}

\begin{proof}
As we have seen already above $p_{12}: W' \longrightarrow S$ has a {\it section} $s: S \longrightarrow W'$ defined by $s(\eta ,\xi)=\stackrel{\sim}{\eta}$ with $\pi(\stackrel{\sim}{\eta})=\eta.$ Since the curve $C_{\eta\xi}$ is rational over $K'_1(\eta ,\xi)$ we have $K'_1(\eta ,\xi)(\omega)\simeq K'_1(\eta ,\xi)(C_{\eta \xi})=K'_1(\eta ,\xi)(\tau)$  with $\tau$ transcendental over $K'_1(\eta ,\xi)$. Therefore we get from (\ref{521})

\begin{equation}\label{523}
K'_1(W')\simeq K'_1(\eta ,\xi ,\omega) \subset K'_1(\eta_1 ,\dots,\eta_{2d-1} ,t,\dots, t_{r-2d+1},\tau )
\end{equation}

\noindent and since the field on the right hand side of (\ref{523}) is purely transcendental over $K'_1$ we have that $W'$ is unirational over $K'_1.$ Finally $p_3$ is surjective, for putting $p_3(\omega )=\alpha$ then $\alpha$ is generic on $\P^r$ over $K'_1$ by corollary \ref{cor8}, hence $\omega $ is generic on $W$ over $K'_1$.
\end{proof}

\bigskip
Finally since $p_3$ is surjective we get the following corollary which concludes the proof of theorem \ref{Thm3} (with $K^{\#}=K'_1$, the field from (\ref{519})).

\begin{corollary}
The variety $W$ is unirational over $K'_1$.
\end{corollary}

\makelastpage

\end{document}